\documentclass{amsart} 
\usepackage{amsmath}
\usepackage[mathscr]{eucal} 
\usepackage{amssymb}
\usepackage{latexsym}
\usepackage{amsthm} 
\theoremstyle{plain}
\newtheorem{theorem}{Theorem\!}[]

\numberwithin{equation}{section}  
\theoremstyle{definition}

\theoremstyle{remark}

\def\Xint#1{\mathchoice
{\XXint\displaystyle\textstyle{#1}}%
{\XXint\textstyle\scriptstyle{#1}}%
{\XXint\scriptstyle\scriptscriptstyle{#1}}%
{\XXint\scriptscriptstyle\scriptscriptstyle{#1}}%
\!\int}
\def\XXint#1#2#3{{\setbox0=\hbox{$#1{#2#3}{\int}$}
\vcenter{\hbox{$#2#3$}}\kern-.5\wd0}}

\def\dashint{\Xint-}

\title[A note on characterization of $H^1$ Sobolev spaces]
{A note on characterization of $H^1$ Sobolev spaces by 
square functions}  
\author{Shuichi Sato} 
 
\begin{document} 
\address{
Kanazawa University,
Kanazawa 920-1192,
Japan}
\email{shuichipm@gmail.com}
\begin{abstract} 
We establish characterization of $H^1$ Sobolev spaces by certain square 
functions, improving previous results.  
  \end{abstract}
  \thanks{2020 {\it Mathematics Subject Classification.\/}
  Primary  42B25; Secondary 46E35. 
  \endgraf
  {\it Key Words and Phrases.} Sobolev spaces, 
Hardy spaces, square functions.   }

\maketitle

\section{ Introduction}  
We use the same notations as in \cite{Sajf}.  
Let $\Phi \in \mathscr M^\alpha$, $\alpha>0$, and 
\begin{align*}\label{ualpha}
U_\alpha(f)(x)&=\left(\int_0^\infty \int_{B_0}  
|f(x-tz) - \Phi_t*f(x-tz)|^2\, dz\, t^{-2\alpha}\frac{dt}{t}
\right)^{1/2}
\\ 
&=\left(\int_0^\infty \int_{B(x,t)} 
| f(z) - \Phi_t* f(z)|^2\, dz\, t^{-2\alpha-n}\frac{dt}{t}
\right)^{1/2},     
\end{align*}  
where $B(x,t)=\{w\in \Bbb R^n: |x-w|<t\}$, $B_0=B(0,1)$. 
We shall prove the following result, which improves \cite[Theorem 1.1]{Sajf}. 

\begin{theorem} \label{thm1}
 Suppose that $n/2<\alpha<n$ and  $\Phi\in \mathscr M^\alpha$. 
Let $f \in  H^1(\Bbb R^n)$. 
Then the following two statements are equivalent:  
\begin{enumerate}  
\item $f\in W^\alpha_{H^1}(\Bbb R^n)$,  
\item $U_\alpha(f) \in L^1(\Bbb R^n)$.  
\end{enumerate} 
Further, we have $\|f\|_{W^\alpha_{H^1}} \simeq 
\|f\|_{H^1} +\|U_\alpha(f)\|_1$.
\end{theorem} 
A version of Theorem \ref{thm1} is shown in \cite{Sajf} by assuming  an additional condition on the size of the Fourier transform of 
$\Phi$ (see \cite[Theorem 1.1]{Sajf}).   
Theorem \ref{thm1} improves \cite[Theorem 1.1]{Sajf} by removing 
the assumption on the Fourier transform of $\Phi$. 
\par 
 Let $\Phi\in \mathscr M^1$.   
Define $\Phi^{(j)}$ as $\Phi^{(1)}(x)=\Phi(x)$ and  
\begin{equation*} 
\Phi^{(j)}(x)=\overbrace{\Phi*\dots *\Phi}^{j}(x), \quad j\geq 2. 
\end{equation*}   
Let $k$ be a positive integer.  Put 
\begin{equation*} 
K^{(k)}(x)=-\sum_{j=1}^k(-1)^j \binom{k}{j}\Phi^{(j)}(x) 
\end{equation*}  
and   
 \begin{equation*} 
 \widetilde{E}_\alpha^{(k)} (f)(x)
=\left(\int_0^\infty \int_{B(x,t)} 
| f(z) - K^{(k)}_t* f(z)|^2\, dz\, t^{-2\alpha-n}\frac{dt}{t}
\right)^{1/2}.   
 \end{equation*} 
Then, since  
$K^{(k)}\in \mathscr M^{\alpha}$ for $\alpha\in (0, 2k)$ 
(see \cite[\S 4]{Sastud}),  by Theorem $\ref{thm1}$  
we have the following.  
\begin{theorem}\label{thm2} 
Suppose that $n/2<\alpha<n$ and $\alpha<\min(2k,n)$. 
Let $f \in  H^1(\Bbb R^n)$. Then,  
$f\in W^\alpha_{H^1}(\Bbb R^n)$ if and only if 
$\widetilde{E}_\alpha^{(k)} (f) \in L^1(\Bbb R^n);$ also 
\begin{equation*}
\|f\|_{W^\alpha_{H^1}} \simeq 
\|f\|_{H^1} +\|\widetilde{E}_\alpha^{(k)} (f)\|_1.  
\end{equation*}  
\end{theorem}  
When $\Phi=|B(0,1)|^{-1}\chi_{B(0,1)}$, $\widetilde{E}_\alpha^{(k)}$ 
is related to 
the operation of repeating the averaging over the ball;  
we observe that if  $k=2$,   
\begin{multline*}
 \widetilde{E}_\alpha^{(2)} (f)(x) 
\\ 
=\left(\int_0^\infty \int\limits_{B(x,t)} \left|f(z)-
 2\dashint_{B(z,t)} f(y)\, dy +
\dashint_{B(z,t)}(f)_{B(y,t)}\, dy\right|^2 
\frac{dz\, dt}{t^{1+2\alpha+n}}\right)^{1/2}, 
\end{multline*} 
where $(f)_{B(y,t)}=\dashint_{B(y,t)} f$.

\section{Proof of Theorem $\ref{thm1}$} 

Let 
\begin{equation*} 
I_{1,1}''(x,y,t)=\int\limits_{z\in B_0, 2|y|/t<|x/t -z|<6} 
\left|(L_\alpha*\Phi)_t(x-y-tz)-(L_\alpha*\Phi)_t(x-tz) \right|^2 \, dz  
\end{equation*} 
be as in \cite[\S 2]{Sajf}, 
where $L_\alpha(x)=\tau(\alpha)|x|^{\alpha-n}$ with  
\begin{equation*}
\tau(\alpha)=\frac{\Gamma\left(n/2-\alpha/2\right)}
{\pi^{n/2}2^\alpha \Gamma\left(\alpha/2\right)}. 
\end{equation*} 
In handling $I_{1,1}''$,  the paper \cite{Sajf} uses 
 assumed size conditions of the Fourier 
transform of $\Phi$. 
In this note we apply different methods and  
under the conditions of Theorem \ref{thm1} we prove 
\begin{equation}\label{e2.1} 
I_{1,1}''(x,y,t)\leq Ct^{-2n}\left(\frac{|y|}{t}\right)^{2\epsilon(n)},  
\quad t>|x|/7>2|y|/7,   
\end{equation}
where $\epsilon(n)=1$ if $n\geq 2$ and $\epsilon(1)=\alpha$.  
Using this modification in the proof of \cite[Lemma 2.1]{Sajf}, 
we can prove \cite[Lemma 2.1]{Sajf} under the conditions that 
$n/2<\alpha<n$ and  $\Phi\in \mathscr M^\alpha$.  Arguing as in the 
proof of \cite[Theorem 1.1]{Sajf} and using 
this improvement on \cite[Lemma 2.1]{Sajf}, we can prove Theorem \ref{thm1}. 
\par 
Thus, to prove Theorem \ref{thm1} it remains to prove \eqref{e2.1}. 
We first note that 
\begin{align*} 
&(L_\alpha*\Phi)_t(x-y-tz)- (L_\alpha*\Phi)_t(x-tz)
\\ 
&=t^{-n}\int_{\Bbb R^n} 
\left(L_\alpha(t^{-1}x-t^{-1}y-z-w)- 
L_\alpha(t^{-1}x-z-w)\right)\Phi(w)\, dw =: J.  
\end{align*} 
Since $\Phi$ is compactly supported, $|x|/t\leq 7$ and $|z|\leq 1$, we have 
\begin{equation*}
J\leq t^{-n}\|\Phi\|_\infty \int_{|w|\leq K} 
\left|L_\alpha(w-t^{-1}y)- L_\alpha(w)\right|\, dw
\end{equation*} 
for some positive constant $K$. Thus 
\begin{multline*}
J\leq Ct^{-n} \int_{|w|\leq 2t^{-1}|y|} 
\left|L_\alpha(w-t^{-1}y)- L_\alpha(w)\right|\, dw 
\\ 
+ Ct^{-n} \int_{2t^{-1}|y|\leq |w|\leq K} 
\left|L_\alpha(w-t^{-1}y)- L_\alpha(w)\right|\, dw 
=: I+II.
\end{multline*} 
We estimate $I$ and $II$ as follows. 
\begin{align*}
I &\leq Ct^{-n} \int_{|w|\leq 2t^{-1}|y|} 
\left|L_\alpha(w-t^{-1}y)\right|\, dw 
 +Ct^{-n} \int_{|w|\leq 2t^{-1}|y|} \left|L_\alpha(w)\right|\, dw 
\\ 
&\leq 2Ct^{-n} \int_{|w|\leq 3t^{-1}|y|} \left|L_\alpha(w)\right|\, dw 
\leq C_\alpha t^{-n} \int_{|w|\leq 3t^{-1}|y|} |w|^{\alpha-n}\, dw
\\ 
&\leq C_\alpha t^{-n} (t^{-1}|y|)^\alpha 
\leq C_\alpha t^{-n} (t^{-1}|y|)^{\epsilon(n)}, 
\end{align*} 
where the last inequality is shown by noting that $\epsilon(n)\leq \alpha$. 
By applying the mean value theorem we see that 
\begin{align*}
II &\leq C_\alpha t^{-n} t^{-1}|y|\int_{2t^{-1}|y|\leq |w|\leq K}
|w|^{\alpha-n-1} \, dw 
\\ 
&\leq C_\alpha t^{-n} t^{-1}|y|\int_{2t^{-1}|y|}^K r^{\alpha-2}\, dr 
\\ 
&\leq C_\alpha t^{-n} (t^{-1}|y|)^{\epsilon(n)}.  
\end{align*} 
Combining the estimates for $I$ and $II$, we have 
$J\leq Ct^{-n} (t^{-1}|y|)^{\epsilon(n)}$, which implies \eqref{e2.1}.

\end{document}